\documentclass[10pt, reqno]{amsart}
\usepackage{amsmath, amsthm, amscd, amsfonts, amssymb, graphicx, color}
\usepackage[bookmarksnumbered, colorlinks, plainpages]{hyperref}
\hypersetup{colorlinks=true,linkcolor=red, anchorcolor=green, citecolor=cyan, urlcolor=red, filecolor=magenta, pdftoolbar=true}
\usepackage{mathrsfs}

\newtheorem{theorem}{Theorem}[section]
\newtheorem{lemma}[theorem]{Lemma}

\newtheorem{corollary}[theorem]{Corollary}
\theoremstyle{definition}

\theoremstyle{remark}

\numberwithin{equation}{section}

\begin{document}
\setcounter{page}{1}

\title[Sklyanin algebras revisited]{Sklyanin algebras revisited}

\author[Nikolaev]
{Igor V. Nikolaev$^1$}

\address{$^{1}$ Department of Mathematics and Computer Science, St.~John's University, 8000 Utopia Parkway,  
New York,  NY 11439, United States.}
\email{\textcolor[rgb]{0.00,0.00,0.84}{igor.v.nikolaev@gmail.com}}





\begin{abstract}
We construct a functor from the category of elliptic curves to a category of 
noncommutative  tori. Our proof is based on an isomorphism between 
the Sklyanin algebras  and   dense 
sub-algebras of the noncommutative tori.  
\end{abstract}

\maketitle

\section{Introduction}
The aim of our note is an elementary proof of the fundamental relation between elliptic curves
and noncommutative tori \cite[Section 1.3]{N}. Such a link  was studied by 
[Sklyanin 1982]   \cite[Section 3]{Skl1},  [Connes \& Dubois-Violette 2002] \cite{CoDu1}
and  [Polishchuk \& Schwarz  2003]  \cite{PoSch1},  but evaded a rigorous proof. 
The key concept of our approach is the Sklyanin algebra $S(\alpha,\beta,\gamma)$.  
By such  one understands  a free   $\mathbf{C}$-algebra   on  four  generators  
$x_i$   satisfying  the following relations: 
\begin{equation}
\left\{
\begin{array}{ccc}
x_1x_2-x_2x_1 &=& \alpha(x_3x_4+x_4x_3),\\
x_1x_2+x_2x_1 &=& x_3x_4-x_4x_3,\\
x_1x_3-x_3x_1 &=& \beta(x_4x_2+x_2x_4),\\
x_1x_3+x_3x_1 &=& x_4x_2-x_2x_4,\\
x_1x_4-x_4x_1 &=& \gamma(x_2x_3+x_3x_2),\\ 
x_1x_4+x_4x_1 &=& x_2x_3-x_3x_2,
\end{array}
\right.
\end{equation}
where $\alpha, \beta$ and $\gamma$ are complex numbers, such that $\alpha+\beta+\gamma+\alpha\beta\gamma=0$
[Smith \& Stafford  1992]   \cite[p. 260]{SmiSta1}.  
 Modulo a two-sided ideal generated by the  central  elements
$\Omega_1 = x_1^2+x_2^2+x_3^2+x_4^2$ and 
$\Omega_2 = x_2^2+{1+\beta\over 1-\gamma}x_3^2+{1-\beta\over 1+\alpha}x_4^2$, 
the  $S(\alpha,\beta,\gamma)$ is the
twisted homogeneous coordinate ring  of an elliptic curve $\mathscr{E}\subset \mathbf{C}P^3$ 
 given in the  Jacobi form, i.e. as an intersection of two quadric surfaces:
\begin{equation}
\left\{
\begin{array}{ccc}
u^2+v^2+w^2+z^2 &=&  0,\\
{1-\alpha\over 1+\beta}v^2+{1+\alpha\over 1-\gamma}w^2+z^2  &=&  0. 
\end{array}
\right.
\end{equation}
In other words, the   $S(\alpha,\beta,\gamma)$ satisfies the fundamental isomorphism
\linebreak
$\hbox{{\bf Mod}}~(S(\alpha,\beta,\gamma))/\hbox{{\bf Tors}}\cong \hbox{{\bf Coh}}~(\mathscr{E})$,
 where {\bf Coh} is  the category of quasi-coherent sheaves on $\mathscr{E}$, 
  {\bf Mod}  the category of graded left modules over the graded ring $S(\alpha,\beta,\gamma)$
 and  {\bf Tors}  the full sub-category of {\bf Mod} consisting of the
torsion modules  [Serre 1955] \cite{Ser1}.  
The algebra $S(\alpha,\beta,\gamma)$ depends on an automorphism
$\sigma: \mathscr{E}\to \mathscr{E}$  [Stafford \& van ~den ~Bergh  2001]  \cite[p. 173]{StaVdb1}.
We assume that $\sigma^4=1$;   in this case  $\beta=1$ and $\gamma=-1$. 
In what follows, we focus on  the Sklyanin algebra $S(\alpha, 1, -1)$.

\medskip
The noncommutative torus $\mathcal{A}_{\theta}$ can be defined as follows. 
Let $S^1$ be the unit circle and  denote by $L^2(S^1)$ the Hilbert space of the square integrable
complex valued functions on $S^1$.  Fix a real number $\theta\in [0,1)$. 
For every $f(e^{2\pi it})\in L^2(S^1)$ we  consider two bounded
linear operators $u$ and $v$  which act by the formulas
$uf(e^{2\pi it})  = f(e^{2\pi i(t-\theta)})$ and 
$vf(e^{2\pi it})  = e^{2\pi it}f(e^{2\pi it})$.
It is verified directly that:  
\begin{equation}\label{eq6}
\left\{
\begin{array}{cc}
vu  &= e^{2\pi i\theta}uv,\\
uu^* &= u^*u = e,\\
vv^* &= v^*v = e,
\end{array}
\right.
\end{equation}
where $u^*$ and $v^*$ are the adjoint operators
and $e$ is the unit operator. The  noncommutative torus 
is defined as  a $C^*$-algebra  $\mathcal{A}_{\theta}$ 
generated by $u$ and $v$ [Wegge-Olsen 1993]  \cite[Section 12.3]{W1}.
By  $\mathcal{A}_{\theta}^0$ we understand 
a dense self-adjoint quotient algebra  of $\mathbf{C}\langle u,u^*, v, v^*\rangle$
modulo a two-sided ideal generated by the relations (\ref{eq6}).
The completion of $\mathcal{A}_{\theta}^0$ in the operator norm on the Hilbert space 
$L^2(S^1)$ is isomorphic to $\mathcal{A}_{\theta}$ [Wegge-Olsen 1993]  \cite[pp. 204-205]{W1}.

\medskip
For brevity,   let  $S(\alpha):=S(\alpha, 1, -1)$ and define a $\ast$-involution on   $S(\alpha)$
by the formulas $x_1\mapsto x_2$ and $x_3\mapsto x_4$. 
Consider a two-sided ideal $I_{\mu}\subset S(\alpha)$
generated by  the relations  $x_1 x_2=x_3 x_4=\mu^{-1}e$, where $e$ is the
unit  of the Sklyanin algebra $S(\alpha)$ and $\mu\in (0,\infty)$ is a real number. 
The symbol  $\cong_\mathbf{C}$ means  an isomorphism of the $\mathbf{C}$-algebras. 
Our main results can be formulated as follows. 
\begin{theorem}\label{thm1}
The formulas $x_1\mapsto u,  x_2\mapsto u^*, x_3\mapsto v, x_4\mapsto v^*$
define a  $\ast$-isomorphism
$S(\alpha) / I_{\mu} \cong_\mathbf{C} ~\mathcal{A}_{\theta}^0$,
where $\theta\in\mathbf{R}/\mathbf{Z},~\mu\in (0,\infty)$ and  $\alpha\in\mathbf{C}$
 depends on $\theta$ and $\mu$. 
\end{theorem}
Denote by  $\{\mathscr{E}\}$  a category of all elliptic curves and 
by  $\{\mathcal{A}_{\theta}\}$ a category of all noncommutative tori. 
Let $F: \{\mathscr{E}\}\to \{\mathcal{A}_{\theta}\}$  be a map defined 
by the isomorphism $S(\alpha) / I_{\mu} \cong_\mathbf{C} ~\mathcal{A}_{\theta}^0$. 
\begin{corollary}\label{cr1}
(i)  $\{\mathscr{E}\}= \{\mathcal{A}_{\theta}\}\times (0,\infty)$ is a trivial 
fiber bundle with projection map  $F: \{\mathscr{E}\}\to \{\mathcal{A}_{\theta}\}$; 
(ii) $F$ is a covariant functor, which maps isomorphic elliptic curves 
$\mathscr{E}, \mathscr{E}'\in  \{\mathscr{E}\}$ to the Morita equivalent 
noncommutative tori $F(\mathscr{E}), F(\mathscr{E}')\in  \{\mathcal{A}_{\theta}\}$.
\end{corollary}

\section{Proofs}
\subsection{Proof of theorem \ref{thm1}}
 We shall split the proof in a  series of lemmas.
\begin{lemma}\label{lm1}
The system of equations (\ref{eq6}) is equivalent to  the following system 
of six quadratic relations: 
\begin{equation}\label{eq8}
\left\{
\begin{array}{ccc}
vu  &=&  e^{2\pi i\theta}uv,\\
uv^* &=&  e^{2\pi i\theta}v^*u,\\
u^*v &=&  e^{2\pi i\theta}vu^*,\\
v^*u^* &=&  e^{2\pi i\theta}u^*v^*,\\
uu^* &=&  u^*u=e,\\
vv^* &= &  v^*v=e.  
\end{array}
\right.
\end{equation}
\end{lemma}
\begin{proof}
Indeed, the first  and the two last equations of (\ref{eq8}) follow immediately from equations
(\ref{eq6}).  We shall proceed stepwise for the rest  of (\ref{eq8}). 

\smallskip
(i) Let us prove that equations  (\ref{eq6}) imply  $uv^*=e^{2\pi i\theta}v^*u$. 
It follows from $uu^*=e$ and $vv^*=e$ that $uu^*vv^*=e$.  Since $uu^*=u^*u$ we can bring  the last 
equation to  the form  $u^*uvv^*=e$ and multiply the  both sides by the constant $e^{2\pi i\theta}$;
thus one gets the equation $u^*(e^{2\pi i\theta}uv)v^*=e^{2\pi i\theta}$.  
But $e^{2\pi i\theta}uv=vu$ and our main equation takes the form $u^*vuv^*= e^{2\pi i\theta}$. 

We can multiply on the left the  both sides of the equation by the element $u$
and thus get the equation $uu^*vuv^*= e^{2\pi i\theta}u$; since $uu^*=e$ 
one  arrives at the equation  $vuv^*= e^{2\pi i\theta}u$.

Again one can multiply on the left the both sides  by the element 
$v^*$ and thus get the equation  $v^*vuv^*= e^{2\pi i\theta}v^*u$; since $v^*v=e$
one gets the required identity  $uv^*= e^{2\pi i\theta}v^*u$.

\smallskip
(ii) Let us prove that equations (\ref{eq6}) imply  $u^*v=e^{2\pi i\theta}vu^*$. 
As in the case (i),  it follows from the equations $uu^*=e$ and $vv^*=e$ that $vv^*uu^*=e$.  Since $vv^*=v^*v$ 
we can bring  the last 
equation to  the form $v^*vuu^*=e$  and multiply the  both sides by the constant $e^{-2\pi i\theta}$;
thus one gets the equation $v^*(e^{-2\pi i\theta}vu)u^*=e^{-2\pi i\theta}$.  
But $e^{-2\pi i\theta}vu=uv$ and our main equation takes the form $v^* uv u^*= e^{-2\pi i\theta}$. 

We can multiply on the left the  both sides of the equation by the element $v$
and thus get the equation $v v^* u v u^*= e^{-2\pi i\theta} v$; since $v v^* =e$ 
one  arrives at the equation  $ u v u^*= e^{-2\pi i\theta} v$.

Again one can multiply on the left the both sides  by the element 
$u^*$ and thus get the equation  $u^* u v u^*= e^{-2\pi i\theta} u^* v$; since $u^* u=e$
one gets the  equation  $v u^*= e^{-2\pi i\theta} u^* v$.  Multiplying both sides
by constant $e^{2\pi i\theta}$ we obtain the required identity  $u^*v=e^{2\pi i\theta} v u^*$.

\smallskip
(iii) Let us prove that equations (\ref{eq6}) imply  $v^* u^* = e^{2\pi i\theta}u^* v^*$. 
Indeed, it was proved in  (i) that $uv^*= e^{2\pi i\theta}v^*u$;  we shall
multiply on the right this equation by the equation $u^*u=e$. Thus one arrives
at  the equation $u v^* u^* u= e^{2\pi i\theta} v^* u $. 

Notice that in the last equation one can cancel $u$ on the right thus bringing 
it to the simpler form  $u v^* u^* = e^{2\pi i\theta} v^*  $.

We shall multiply on the left both sides of the above equation by the element $u^*$;
one gets therefore $u^* u v^* u^* = e^{2\pi i\theta} u^* v^*  $.  But $u^* u = e$
and the left hand side simplifies giving  the required identity
  $v^* u^* = e^{2\pi i\theta} u^* v^*  $.
\end{proof}

\begin{lemma}\label{lm2}
Each  Sklyanin  algebra $S(\alpha)$ is isomorphic to
a free algebra 
\linebreak
$\mathbf{C}\langle u,u^*, v, v^*\rangle$  modulo an ideal generated  by  six  
skew-symmetric quadratic  relations:
\begin{equation}\label{eq9}
\left\{
\begin{array}{ccc}
v u &=& \mu e^{2\pi i\theta} u v,\\
u v^*  &=& {1\over\mu} e^{2\pi i\theta}   v^* u,\\
u^* v  &=& \mu e^{2\pi i\theta} v u^* ,\\
v^* u^* &=& {1\over \mu} e^{2\pi i\theta} u^* v^*,\\
 u  u^*  &=&  u^* u ,\\
v v^*    &=& v^* v,
\end{array}
\right.
\end{equation}
where $\theta\in [0, 1)$ and $\mu\in (0,\infty)$. 
\end{lemma}
\begin{proof}
(i)  Since for  the Sklyanin algebra $S(\alpha)$ we have $\sigma^4=1$,
such an algebra  is isomorphic to the  algebra  $\mathbf{C}\langle u,u^*,v,v^*\rangle$
modulo  six skew-symmetric relations: 
\begin{equation}\label{eq10}
\left\{
\begin{array}{ccc}
vu &=& q_{13} uv,\\
v^*u^* &=&  q_{24}u^*v^*,\\
v^*u &=&  q_{14}uv^*,\\
vu^* &=&  q_{23}u^*v,\\
u^*u &=&  q_{12}uu^*,\\
v^*v &=&  q_{34}vv^*,
\end{array}
\right.
\end{equation}
where $q_{ij}\in \mathbf{C}\setminus\{0\}$,  see  [Feigin \& Odesskii  1989]  \cite[Remark 1]{FeOd1}
and   [Feigin \& Odesskii  1993]  \cite[Section 2]{FeOd2} for the details.

\bigskip
(ii) It is verified directly,  that relations (\ref{eq10})  are invariant of the involution, 
if and only if,  the following restrictions on the constants $q_{ij}$ hold:
\begin{equation}\label{eq11}
\left\{
\begin{array}{ccc}
q_{13} &=&  (\bar q_{24})^{-1},\\
q_{24} &=&  (\bar q_{13})^{-1},\\
q_{14} &= & (\bar q_{23})^{-1},\\
q_{23} &= & (\bar q_{14})^{-1},\\
q_{12} &= & \bar q_{12},\\
q_{34} &= & \bar q_{34},
\end{array}
\right.
\end{equation}
where $\bar q_{ij}$ means the complex conjugate of $q_{ij}\in \mathbf{C}\setminus\{0\}$.

\bigskip
(iii)
Consider a one-parameter family  $S(q_{13})$ of the Sklyanin $\ast$-algebras
defined by relations (\ref{eq10})  with  restrictions (\ref{eq11}),   where 
$\bar q_{14}=q_{13},$ and $q_{12} =  q_{34}=1$.
It is not hard to see,  that the $\ast$-algebras  $S(q_{13})$ 
are pairwise non-isomorphic for different values of complex  parameter $q_{13}$;
therefore  the family $S(q_{13})$ is   a
normal form of  the Sklyanin $\ast$-algebra $S(\alpha)$. 
It remains to notice, that one can write  complex parameter $q_{13}$
in the polar form $q_{13}=\mu e^{2\pi i\theta}$,  where $\theta=Arg~(q_{13})$
and $\mu=|q_{13}|$.   Lemma \ref{lm2} follows.
\end{proof}

\begin{lemma}\label{lm4}
The system of relations (\ref{eq8}) is equivalent to relations (\ref{eq9})
plus the scaled unit relation  $uu^*=vv^*={1\over\mu}e$.
\end{lemma}
\begin{proof}
(i) Using the last two relations,  one can bring (\ref{eq8}) to the form:
\begin{equation}\label{eq13}
\left\{
\begin{array}{ccc}
vuv^* &=&  e^{2\pi i\theta}u,\\
v^* &= & e^{2\pi i\theta}u^*v^*u,\\
v^*uv &=&  e^{-2\pi i\theta}u,\\
u^* &=&  e^{-2\pi i\theta}v^*u^*v,\\
uu^*   &=&   u^*u   =e,\\
 vv^*    &=&  v^*v  =e.
\end{array}
\right.
\end{equation}

\bigskip
(ii)  If one adjoins relation $uu^*=vv^*={1\over\mu}e$ to 
(\ref{eq9}),  then $\mu$ can be eliminated  from
the first four relations and the system of relations (\ref{eq9}) takes the form: 
\begin{equation}\label{eq14}
\left\{
\begin{array}{ccc}
vuv^* &=&  e^{2\pi i\theta}u,\\
v^* &=&  e^{2\pi i\theta}u^*v^*u,\\
v^*uv &=&  e^{-2\pi i\theta}u,\\
u^* &= & e^{-2\pi i\theta}v^*u^*v,\\
uu^*   &=&   u^*u ={1\over\mu}e,\\
 vv^*    &=&  v^*v={1\over\mu}e.
\end{array}
\right.
\end{equation}

\bigskip
(iii)  Comparing systems (\ref{eq13}) and (\ref{eq14}),   one concludes
that they  coincide up to a scaled unit $e'={1\over\mu} e$.
Lemma \ref{lm4} follows.  
\end{proof}

\bigskip\noindent
Theorem \ref{thm1} follows from lemma \ref{lm4},  since the ideal $I_{\mu}$
of the Sklyanin algebra $S(\alpha)$ is generated  by relations 
 $uu^*=vv^*={1\over\mu}e$.

\subsection{Proof of corollary \ref{cr1}}
(i) The fact that $\{\mathscr{E}\}$ fibers over $\{\mathcal{A}_{\theta}\}$ 
follows from part (iii) of proof of lemma \ref{lm2}. Indeed, the set 
 $\{\mathscr{E}\}$ is parametrized by the complex variable $q_{13}=\mu e^{2\pi i\theta}$. 
Since $\theta\in\mathbf{R}/\mathbf{Z}$ and $\mu\in (0,\infty)$, we conclude that
$\{\mathscr{E}\}= \{\mathcal{A}_{\theta}\}\times (0,\infty)$. 

\smallskip
(ii) If $\mathscr{E},\mathscr{E}'\in\{\mathscr{E}\}$ are isomorphic elliptic curves,
then $\hbox{{\bf Coh}}~(\mathscr{E})\cong \hbox{{\bf Coh}}~(\mathscr{E}')$ are
equivalent categories. The Serre isomorphism  $\hbox{{\bf Coh}}~(\mathscr{E})\cong 
\hbox{{\bf Mod}}~(S(\alpha))/\hbox{{\bf Tors}}$ says that 
$\hbox{{\bf Mod}}~(S(\alpha))\cong \hbox{{\bf Mod}}~(S'(\alpha))$ are equivalent categories. 
But any equivalence in the category of  modules corresponds to the Morita equivalence
of the underlying Sklyanin algebras $S(\alpha)$. In view of theorem \ref{thm1},
the same is true for the algebras  $\mathcal{A}_{\theta}^0$ and their operator norm closures $\mathcal{A}_{\theta}$.

\bibliographystyle{amsplain}


\end{document}